\documentclass[11pt,reqno]{amsart}
\usepackage[a4paper,margin=27mm]{geometry}
\usepackage{amsmath,amssymb,mathtools}
\usepackage[T1]{fontenc}
\usepackage{lmodern}
\usepackage[expansion=false]{microtype}
\usepackage{hyperref}
\theoremstyle{plain}
\newtheorem{theorem}{Theorem}
\newtheorem{lemma}[theorem]{Lemma}
\newtheorem{proposition}[theorem]{Proposition}
\newtheorem{corollary}[theorem]{Corollary}
\theoremstyle{remark}\newtheorem{remark}[theorem]{Remark}

\newcommand{\E}{\mathbb E}
\newcommand{\Prb}{\mathbb P}
\newcommand{\R}{\mathbb R}
\newcommand{\K}{\mathbb K}

\title[Very tight embeddings below exponent one]{Very tight embeddings of subspaces of \(L_p\) for \(0<p<1\)}

\author{Gang Cai}
\address[G. Cai]{School of Mathematical Sciences, Chongqing Normal University, Chongqing 401331, China}
\email{caigang-aaaa@163.com}

\author[Longfa Sun]{Longfa Sun$^\dag$}
\address[L. Sun]{Hebei Key Laboratory of Physics and Energy Technology, School of Mathematics and Physics, North China Electric Power University, Baoding 071003, China}
\email{sun.longfa@ncepu.edu.cn}

\thanks{$^\dag$ Corresponding author.\\
This work is supported by the National Natural Science Foundation of China (Grant No. 12571141), the
Chongqing Talent Support Program (Grant No. cstc2024ycjh-bgzxm0121), and the Fundamental Research Funds for the Central Universities (Grant No. 2025MS178).}

\date{}
\makeatletter
\@ifundefined{subjclassname@2020}{%
  \@namedef{subjclassname@2020}{\textup{2020} Mathematics Subject Classification}%
}{}
\makeatother
\subjclass[2020]{Primary 46B07; Secondary 46B09, 46E30}
\keywords{Finite-dimensional subspaces, quasi-Banach spaces,
linear embeddings, Rademacher type, change of density}
\begin{document}
\begin{abstract}
For \(0<p<1\) and \(\mu>1\), every \(d\)-dimensional subspace
of \(L_p\) embeds linearly into \(\ell_p^{\lfloor\mu d\rfloor}\)
with distortion at most \(C(p,\mu)(\log(d+2))^{1/p}\).
If its Rademacher type \(q\) constant is at most \(T\), where
\(p<q\leq2\), the distortion is at most \(C(p,q,\mu)T\),
independently of \(d\). We also give explicit bounds as
\(\mu\downarrow1\). This answers the question of Johnson and
Schechtman concerning the extension of their very tight embedding
theorem below exponent one. The proof combines compression at a
lower exponent with a change of density that preserves the number
of coordinates, and applies over both scalar fields.
\end{abstract}
\maketitle

\section{Introduction}
The discretization of finite-dimensional subspaces of \(L_p\)
asks how many coordinates are needed to preserve their geometry
with controlled distortion. The least possible number is the
dimension of the source. A very tight embedding approaches this
lower bound: for every fixed \(\mu>1\), a \(d\)-dimensional
source is represented in \(\lfloor\mu d\rfloor\) coordinates.
The distortion may depend on \(\mu\). It need not approach one
as \(\mu\downarrow1\).

Discretization results for subspaces of \(L_p\) were developed in
\cite{BLM89,Sch87,Tal90,Tal95}, and in the range \(0<p<1\) in
\cite{Pena98,SZ01,Zv00}. An account of the theory is given in
\cite{JS01}. Johnson and Schechtman \cite{JS03} obtained very tight
embeddings for \(1\leq p<2\). Their Theorem~2 gives
dimension-independent distortion under an improved type assumption,
and distortion of order \((\log(d+2))^{1/p}\) for arbitrary
subspaces. They asked whether these conclusions remain valid for
\(0<p<1\) in \cite[Remark~2, p.~849]{JS03}. We prove that they
do, and obtain quantitative dependence on the type constant and
the excess of target dimension.

We work over \(\K=\R\) or \(\mathbb C\), with dimensions taken
over \(\K\). For \(s>0\), the spaces \(L_s(\Omega,\lambda;\K)\)
and \(\ell_s^n=\K^n\) carry the usual functionals
\[
 \|f\|_s=\left(\int_\Omega |f|^s\,d\lambda\right)^{1/s},
 \qquad \|x\|_s=\left(\sum_{j=1}^n|x_j|^s\right)^{1/s}.
\]
Subspaces inherit these functionals. A linear embedding
\(A:X\to Z\) satisfies
\(c\|x\|_X\leq\|Ax\|_Z\leq C\|x\|_X\) for some
\(0<c\leq C<\infty\). Its distortion is
\[
 \operatorname{dist}(A)=
 \frac{\sup_{x\ne0}\|Ax\|_Z/\|x\|_X}
      {\inf_{x\ne0}\|Ax\|_Z/\|x\|_X}.
\]
A \(D\)-embedding has distortion at most \(D\). An isometric
embedding preserves the functional exactly. Constants depend only
on the displayed parameters and may increase between occurrences.
We use the usual floor and ceiling notation, and \(\log\) denotes
the natural logarithm.

For \(0<q\leq2\), the Rademacher type \(q\) constant \(T_q(X)\)
is the least \(C\) for which
\begin{equation}\label{eq:type-definition}
 \left(\E\left\|\sum_{i=1}^n\varepsilon_i x_i\right\|_X^2\right)^{1/2}
 \leq C\left(\sum_{i=1}^n\|x_i\|_X^q\right)^{1/q}
\end{equation}
holds for every finite sequence in \(X\), where the
\(\varepsilon_i\) are independent uniform signs. We set
\(T_q(X)=\infty\) if no finite constant works. Thus a bound
\(T_q(X)\leq T\) uses the same \(T\) for all finite sequences.
In the following theorem, \(T\) is also uniform over the source
spaces and their dimensions.

\begin{theorem}\label{thm:targets}
Let \(0<p<1\), \(\mu>1\), let \(d\geq1\) be an integer, and put
\(m=\lfloor\mu d\rfloor\). The following assertions hold over
both scalar fields.
\begin{enumerate}
\item[(a)] If \(p<q\leq2\), \(1\leq T<\infty\), and
\(T_q(X)\leq T\), then every
\(d\)-dimensional subspace \(X\subset L_p\) admits a linear
embedding into \(\ell_p^m\) with distortion at most
\(C(p,q,\mu)T\).
\item[(b)] Every \(d\)-dimensional subspace \(X\subset L_p\)
admits a linear embedding into \(\ell_p^m\) with distortion at
most \(C(p,\mu)(\log(d+2))^{1/p}\).
\end{enumerate}
For \(1<\mu\leq2\), the respective distortion bounds can be
chosen as
\begin{equation}\label{eq:main-quantitative}
 C(p,q)T\exp\!\left(\frac{C(p)}{\mu-1}\right),
 \qquad
 \exp\!\left(\frac{C(p)}{\mu-1}\right)(\log(d+2))^{1/p}.
\end{equation}
\end{theorem}

The target dimension is prescribed for every \(\mu>1\).
A discretization in \(C_0d\) coordinates for one fixed \(C_0\)
does not suffice. The random operator used below is the
deterministic-weight construction of Friedland and Gu\'edon
\cite{FG11}, with related estimates in \cite{BLV07}.
Their whole-space theorem embeds \(\ell_a^n\) into \(\ell_r^M\)
for \(M>n\), with distortion bounded by
\(C(a,r)^{n/(M-n)}\) in its range of exponents.
In their subsequent work \cite{FG13}, they obtained the stronger
bound
\[
 C(a,r)\left(\frac{\log(1+1/\theta)}{\theta}\right)^{1/r-1/a}
\]
for real whole-space embeddings
\(\ell_a^n\to\ell_r^{(1+\theta)n}\), where
\(0<r<a<2\), \(r\leq1\), and \(0<\theta\leq1\),
with the usual integer convention for the target dimension.
That construction combines a scaled identity on \(n\) coordinates
with a random operator on the additional coordinates.

Here the source is a prescribed subspace \(E\subset\ell_a^N\)
of dimension \(d\), and the required number of target coordinates
is \(m=\lfloor\mu d\rfloor\), which may be smaller than \(N\).
Restricting either whole-space embedding to \(E\) retains a target
dimension determined by \(N\), rather than by \(d\).
In particular, retaining the identity block from \cite{FG13}
already uses \(N\) coordinates. Those conclusions therefore do
not directly yield the subspace compression needed here.
Proposition~\ref{prop:compression} proves this compression when
\(N\leq C_0d\). Our exponential bound in \((\mu-1)^{-1}\)
is not an improvement of the whole-space bound in \cite{FG13}.
The point is to control the target dimension by the dimension of
the prescribed subspace.

We use the mean, concentration, and Fourier small-ball estimates
for this deterministic-weight random operator in the form proved
below. The construction is taken from \cite{FG11}. Our argument
keeps the dependence on \(N/m\) explicit in the upper and
small-ball bounds, and uses the dimension of \(E\) only in the
final lower-bound net. We give the estimates over both scalar
fields so that they apply also when \(m<N\).
The proof controls the random operator above on the ambient space
and below on \(E\). The small-ball estimate has power
\(\sigma m\), where \(\sigma=\dim_{\R}\K\), whereas a net
in \(E\) contributes power \(\sigma d\). Their combination
leaves \(\sigma(m-d)\), allowing any fixed proportional excess
of target dimension. The estimates use the quasi-triangle
inequalities and do not require convex unit balls.

Both parts of Theorem~\ref{thm:targets} use the same initial
construction. Lemma~\ref{lem:preliminary}, with \(a=p/2\),
and Proposition~\ref{prop:compression}, with \(r=3p/8\), give
\[
 X\longrightarrow E\subset\ell_{p/2}^N
 \longrightarrow Y\subset\ell_{3p/8}^m,
 \qquad N\leq C(p)d.
\]
The distortion of this map depends only on \(p,\mu\).
Lemma~\ref{lem:restore} then restores exponent \(p\) without
adding coordinates. In (a), \(p\) is below the available type
exponent. In (b), it equals that exponent. The latter passage
uses a regularization of the atom weights and produces the
logarithmic factor. Applying the density theorem to the image of
the original unit ball retains linear
dependence on both the type constant and the distortion of the
initial embedding, and gives \eqref{eq:main-quantitative}.
The preliminary results, compression estimates, and density
argument are presented in the next three sections, followed by
the proof of the theorem.

\section{Preliminaries}\label{sec:preliminaries}
We retain the scalar-field and dimension conventions from the
introduction and put \(\sigma=\dim_{\R}\K\in\{1,2\}\).
Thus a \(d\)-dimensional space over \(\K\) has real dimension
\(\sigma d\). Fourier transforms and volumes use this real dimension. The scalar field is omitted from the notation \(L_p\) and
\(\ell_p^n\) when its choice is clear. Constants can be chosen
uniformly for these two scalar fields.
A quasi-norm on a vector space \(X\) over \(\K\) is a function
\(\|\cdot\|_X:X\to[0,\infty)\) such that \(\|x\|_X=0\) if and
only if \(x=0\), \(\|\lambda x\|_X=|\lambda|\|x\|_X\), and
\(\|x+y\|_X\leq K_X(\|x\|_X+\|y\|_X)\) for a fixed
\(K_X\geq1\). A complete quasi-normed space is called a
quasi-Banach space. For \(0<p<1\), the usual \(L_p\) and
\(\ell_p\) functionals are quasi-norms satisfying the
\(p\)-triangle inequality.

An embedding is a continuous linear isomorphism onto its image.
For a bounded linear map, we use the operator norm
\[
 \|A\|=\sup_{x\ne0}\frac{\|Ax\|_Z}{\|x\|_X}.
\]
The distortion defined in the introduction can equivalently be written as
\[
 \operatorname{dist}(A)=\|A\|\,\|A^{-1}\|.
\]
For composable embeddings,
\[
 \operatorname{dist}(BA)\leq
 \operatorname{dist}(B)\operatorname{dist}(A),
\]
by submultiplicativity of the operator norms of the maps and their
inverses. Inverse maps are always taken on their images.
For a quasi-normed space \(E\), write
\(B_E=\{x\in E:\|x\|_E\leq1\}\) for its closed unit ball.
For a measurable set \(D\), \(1_D\) denotes its indicator function,
and \(\operatorname{span}_{\K}S\) denotes the finite linear
combinations of elements of \(S\) with coefficients in \(\K\).
An \(a\)-norm is a homogeneous quasi-norm satisfying
\(\|x+y\|^a\leq\|x\|^a+\|y\|^a\).
We use this terminology for \(0<a\leq1\). In this case
\(d(x,y)=\|x-y\|^a\) is a metric, and
\[
 \bigl|\|x\|^a-\|y\|^a\bigr|\leq\|x-y\|^a.
\]
Indeed, apply the \(a\)-triangle inequality to \(x=y+(x-y)\)
and then interchange \(x\) and \(y\). This replacement for the
usual reverse triangle inequality will be used in the random-operator
estimates. A \(\delta\)-net of a set \(S\) is a subset
\(\mathcal N\subset S\) such that for every \(x\in S\) there
is \(y\in\mathcal N\) with \(\|x-y\|\leq\delta\).
Thus its radius is \(\delta^a\) in the metric just defined.

We shall repeatedly normalize embeddings as follows. If
\(A_0:X\to Z\) is an embedding of a nonzero space with distortion
at most \(D\), set \(A=\|A_0^{-1}\|A_0\). Then
\[
 \|x\|_X\leq\|Ax\|_Z\leq D\|x\|_X.
\]
The first inequality follows from
\(\|x\|_X\leq\|A_0^{-1}\|\|A_0x\|_Z\). The second follows
from the definition of distortion. Multiplication by this scalar
changes neither the image dimension nor the distortion.

We use the Rademacher type constants defined in
\eqref{eq:type-definition}. We record two consequences needed below:
the automatic type bound for \(L_q\) when \(q\leq1\), and the
transport of type under a linear embedding.
For \(0<q\leq1\), \(L_q\) has type \(q\) constant at most one, by its
\(q\)-triangle inequality for each choice of signs.

We shall also use the elementary transport inequality
\begin{equation}\label{eq:type-transport}
 T_q(A(X))\leq \operatorname{dist}(A)T_q(X).
\end{equation}
To prove \eqref{eq:type-transport}, let \(y_1,\ldots,y_n\in A(X)\)
and set \(x_i=A^{-1}y_i\). The quasi-norm on \(A(X)\) is inherited
from \(Z\), and linearity gives
\(\sum_i\varepsilon_i y_i=A(\sum_i\varepsilon_i x_i)\).
Consequently,
\begin{align*}
 \left(\E\left\|\sum_{i=1}^n\varepsilon_i y_i\right\|_Z^2\right)^{1/2}
 &\leq \|A\|
   \left(\E\left\|\sum_{i=1}^n\varepsilon_i x_i\right\|_X^2\right)^{1/2}\\
 &\leq \|A\|T_q(X)
   \left(\sum_{i=1}^n\|x_i\|_X^q\right)^{1/q}\\
 &\leq \|A\|\,\|A^{-1}\|T_q(X)
   \left(\sum_{i=1}^n\|y_i\|_Z^q\right)^{1/q}.
\end{align*}
Since the sequence \((y_i)\) is arbitrary, the definition of
\(T_q(A(X))\) gives \eqref{eq:type-transport}.
If \(T_q(X)=\infty\), the assertion is immediate.

For a probability measure \(\rho\) and \(0<q<\infty\), write
\[
 \|h\|_{q;\rho}=\left(\int|h|^q\,d\rho\right)^{1/q},
 \qquad
 \|h\|_{q,\infty;\rho}
 =\sup_{t>0}t\,\rho\{|h|>t\}^{1/q}.
\]
The first is the usual (strong) \(L_q\) functional, and the second
is the weak \(L_q\) quasi-norm. The space \(L_{q,\infty}(\rho)\)
consists of the functions for which the latter is finite. The
function \(t\mapsto\rho\{|h|>t\}\) is the tail distribution function
of \(|h|\) used below. If the measure is omitted from a quasi-norm,
it is the specified ambient measure. In a change of density to
\(g\,d\nu\), the quasi-norm of the transformed function is taken
with respect to that new measure.

We first record the following established results over the real field.
The first gives a number of coordinates proportional to the source
dimension, the second changes integrability after a fixed multiplication,
and the third lowers the exponent while preserving the source geometry.
Only these stated conclusions will be used from the cited results.
\begin{enumerate}
\item \label{input:discretization}\textbf{Preliminary discretization}
\cite[Theorem~13(i), Section~2.1]{JS01}:
let \(0<a<q<2\), \(0<T<\infty\), and \(\epsilon>0\).
If \(X\subset L_a\) is \(d\)-dimensional and \(T_q(X)\leq T\),
then there exists a linear embedding \(J:X\to\ell_a^N\) such that
\[
 \operatorname{dist}(J)\leq1+\epsilon,
 \qquad N\leq C(a,q,T,\epsilon)d.
\]
The case \(0<a<1\) is included in Pe\~na's result \cite{Pena98}.
With the embedding convention from the introduction, the strict inequality \(a<q\) is the hypothesis
that permits the linear bound on \(N\). We use this result only with
\(\epsilon=1\), so no estimate of the dependence on \(\epsilon\)
is needed. The type convention above is the second-moment convention
used in the survey. In Lemma~\ref{lem:preliminary}, the isometric
change in item~\ref{input:stable} places the original space in
\(L_a\), where its automatic type \(p\), with \(a<p\), provides
the required exponent gap.
\item \label{input:density}\textbf{Pisier's density theorem}
\cite[Theorem~11]{JS01}, in the factorization framework of
Maurey \cite{Maurey74} and Pisier \cite{Pisier86}:
we use the implication from the finite-family maximum estimate to
a density with a local integral estimate in
\cite[Theorem~11]{JS01}. Written with source exponent \(r\) and
target exponent \(q\), the needed statement is the following:
if \(0<r<q<\infty\), \(S\subset L_r(\nu)\), \(\nu\) is a probability,
and
\[
 \left\|\max_i|c_i f_i|\right\|_r
       \leq M\left(\sum_i|c_i|^q\right)^{1/q}
 \quad(f_i\in S,\ c_i\in[0,1]),
\]
for every finite family \((f_i)\) and every corresponding choice of
coefficients \((c_i)\), then there exists a probability density \(g\) with respect to \(\nu\),
that is, \(g\geq0\) and \(\int g\,d\nu=1\), such that the map
\(Uf=g^{-1/r}f\), restricted to \(\operatorname{span}_{\R}S\),
is an \(L_r\)-isometry into \(L_r(g\,d\nu)\) and
\(\|Uf\|_{q,\infty;g\,d\nu}\leq C(r,q)M\) for \(f\in S\).
The dependence on \(M\) is linear. For \(M>0\), apply the
theorem to \(M^{-1}S\), whose maximum estimate has constant one,
and then multiply the conclusion by \(M\). The same density
works for \(S\). If \(M=0\), every element of \(S\) is zero
almost everywhere and the assertion is immediate.
Here every \(f\in S\) vanishes \(\nu\)-almost everywhere on
\(\{g=0\}\), and \(Uf\) is defined to be zero there.
This support assertion follows from the equivalent local estimate
in the cited theorem:
\[
 \|1_Ef\|_r\leq C\left(\int_Eg\,d\nu\right)^{1/r-1/q}
 \qquad(f\in S).
\]
Here \(E\) ranges over all measurable sets and
\(C=C(r,q)M\), after increasing the constant depending on \(r,q\).
Choose \(g\) from this local estimate. Taking \(E=\{g=0\}\)
proves the support assertion, which also holds on \(\operatorname{span}_{\R}S\).
For \(E_t=\{|g^{-1/r}f|>t\}\), the same estimate gives
\[
 t^r\int_{E_t}g\,d\nu
 \leq\int_{E_t}|f|^r\,d\nu
 \leq C^r\left(\int_{E_t}g\,d\nu\right)^{1-r/q}.
\]
Write \(b_t=\int_{E_t}g\,d\nu\). If \(b_t>0\), dividing the
last inequality by \(b_t^{1-r/q}\) and taking the \(r\)-th root
gives \(t b_t^{1/q}\leq C\). This also holds when \(b_t=0\).
Taking the supremum over \(t>0\) proves the weak \(L_q\) bound.
For any \(h\in\operatorname{span}_{\R}S\), the support assertion gives
\[
 \int |Uh|^r g\,d\nu
 =\int_{\{g>0\}} |g^{-1/r}h|^r g\,d\nu
 =\int |h|^r\,d\nu.
\]
Thus the change of measure and the multiplier cancel exactly in the
\(L_r\) quasi-norm. Both are fixed for the entire set \(S\).

\item \label{input:stable}\textbf{Isometric embedding via stable random variables:}
if \(0<a<p<2\), a finite-dimensional subspace of \(L_p\) embeds
linearly isometrically into an \(L_a\)-space. It follows from the
classical stable representation \cite{BDK66}, also recalled in \cite{JS01}.
The input from stable probability is the existence of the random
vector with the characteristic function displayed below, together with
finiteness of its one-dimensional moments of order less than \(p\).
The deduction of the isometric embedding from these facts is included
here to specify its normalization and its linearity.
For a real random vector \(\Theta\), its characteristic function is
\(t\mapsto\E e^{i\langle t,\Theta\rangle}\), where
\(\langle t,u\rangle=\sum_jt_ju_j\). A symmetric \(p\)-stable
vector is one whose every real linear functional has characteristic
function \(u\mapsto\exp(-c|u|^p)\) for some \(c\geq0\).
The formula below specifies this constant for each functional.

More explicitly, fix a basis \(f_1,\ldots,f_d\) of \(X\). There exists a
symmetric \(p\)-stable random vector \(\Theta\in\R^d\)
with characteristic function
\[
 \E e^{i\langle t,\Theta\rangle}
   =\exp\left(-\left\|\sum_{j=1}^d t_jf_j\right\|_p^p\right).
\]
If \(Z\) has characteristic function \(e^{-|s|^p}\), put
\(b_{a,p}=(\E|Z|^a)^{1/a}\in(0,\infty)\). The map
\[
 I\left(\sum_j t_jf_j\right)=b_{a,p}^{-1}\langle t,\Theta\rangle
\]
is well defined and real linear, since the basis coefficients are
unique and the same random vector \(\Theta\) is used for every \(x\).
For \(x=\sum_jt_jf_j\) and \(u\in\R\), homogeneity gives
\[
 \E e^{iu\langle t,\Theta\rangle}
 =\exp\left(-\left\|\sum_jut_jf_j\right\|_p^p\right)
 =\exp(-|u|^p\|x\|_p^p).
\]
This is the characteristic function of \(\|x\|_pZ\). Uniqueness
of characteristic functions therefore identifies their distributions.
Since \(a<p\), the moment defining \(b_{a,p}\) is finite and
strictly positive. Consequently \(Ix\in L_a\), and
\[
 \|Ix\|_a^a
 =b_{a,p}^{-a}\E|\langle t,\Theta\rangle|^a
 =b_{a,p}^{-a}\|x\|_p^a\E|Z|^a
 =\|x\|_p^a.
\]
\end{enumerate}
The next lemma gives the preliminary reduction over both scalar
fields. We first lower the exponent isometrically and then discretize
using the automatic type \(p\) of the original space. For complex
spaces, realification is followed by an explicit complex-linear
construction that preserves the required coordinate bound.

\begin{lemma}\label{lem:preliminary}
Let \(0<a<p<1\), and let \(X\subset L_p\) have dimension
\(d\) over \(\K\). There is a \(\K\)-linear embedding
\(J:X\to\ell_a^N\) with distortion at most
\(D_a=2^{1+1/a}\) and \(N\leq C(a,p)d\).
\end{lemma}
\begin{proof}
The zero-dimensional case is immediate, so assume \(d\geq1\).
First suppose that \(\K=\R\).
Since \(0<a<p<1<2\), the stable representation in item~\ref{input:stable}
applies to the given real subspace of \(L_p\). It gives a real-linear isometry
\(I:X\to L_a(\rho;\R)\) on a suitable probability space.
For every finite sequence \((x_i)\) in \(X\), the \(p\)-triangle
inequality in the original \(L_p\)-space yields, for every choice of signs,
\[
 \left\|\sum_i\varepsilon_i Ix_i\right\|_a
 =\left\|\sum_i\varepsilon_i x_i\right\|_X
 \leq\left(\sum_i\|x_i\|_X^p\right)^{1/p}
 =\left(\sum_i\|Ix_i\|_a^p\right)^{1/p}.
\]
Taking the second moment proves \(T_p(IX)\leq1\). We can therefore
apply preliminary discretization to \(IX\subset L_a\), with type
exponent \(p\) and \(\varepsilon=1\), since \(a<p<1<2\). Composing the resulting
embedding with \(I\), and normalizing, gives a real-linear map
\(J_{\R}:X\to\ell_a^N(\R)\), with \(N\leq C(a,p)d\), such that
\begin{equation}\label{eq:preliminary-real}
 \|x\|_X\leq\|J_{\R}x\|_a\leq2\|x\|_X
 \qquad(x\in X).
\end{equation}

Now let \(\K=\mathbb C\), and write \(X_{\R}\) for the same
space with scalars restricted to \(\R\). We describe explicitly an
isometric realization of the realification of the given complex \(L_p\)-space.
Restricting the original measure to the union of the supports of a
finite basis of \(X\), we may assume that \((\Omega,\lambda)\)
is \(\sigma\)-finite: for each basis function \(f\), the sets
\(\{|f|>1/n\}\) have finite measure and cover its support.
Let \(\tau\) be normalized Lebesgue measure on \([0,2\pi)\), and put
\[
 c_p=\left(\int_0^{2\pi}|\cos t|^p\,d\tau(t)\right)^{1/p}>0.
\]
Define the real-linear map
\[
 (Qf)(\omega,t)=c_p^{-1}\operatorname{Re}(e^{it}f(\omega)).
\]
For every \(z\in\mathbb C\), rotational invariance of \(\tau\)
gives
\[
 \int_0^{2\pi}|\operatorname{Re}(e^{it}z)|^p\,d\tau(t)
 =c_p^p|z|^p.
\]
Consequently, Tonelli's theorem gives
\begin{align*}
 \|Qf\|_{L_p(\lambda\times\tau;\R)}^p
 &=c_p^{-p}\int_\Omega\int_0^{2\pi}
       |\operatorname{Re}(e^{it}f(\omega))|^p\,d\tau(t)\,d\lambda(\omega)\\
 &=c_p^{-p}\int_\Omega c_p^p|f(\omega)|^p\,d\lambda(\omega)
 =\|f\|_{L_p(\lambda;\mathbb C)}^p.
\end{align*}
In particular, \(Q\) is injective and realizes \(X_{\R}\)
isometrically as a real \(L_p\) subspace of dimension \(2d\).
Apply the real assertion to \(Q(X_{\R})\), and compose with
\(Q\). This gives a real-linear map
\(J_{\R}:X_{\R}\to\ell_a^N(\R)\) satisfying
\eqref{eq:preliminary-real}, with \(N\leq2C(a,p)d\).
The map \(J_{\R}\) need not commute with multiplication by \(i\).
The following construction enforces this identity while keeping
exactly \(N\) complex coordinates.
Define a single map into \(\ell_a^N(\mathbb C)\) by
\[
 Jx=J_{\R}x-iJ_{\R}(ix).
\]
This map is real linear, and
\[
 J(ix)=J_{\R}(ix)-iJ_{\R}(-x)
      =J_{\R}(ix)+iJ_{\R}x=iJx.
\]
It is therefore complex linear. Since the coordinates of \(J_{\R}x\)
and \(J_{\R}(ix)\) are real, each coordinate satisfies
\[
 |(Jx)_\ell|\geq|(J_{\R}x)_\ell|,
 \qquad
 |(Jx)_\ell|^a
 \leq|(J_{\R}x)_\ell|^a+|(J_{\R}(ix))_\ell|^a.
\]
Summing these inequalities and using \eqref{eq:preliminary-real},
together with \(\|ix\|_X=\|x\|_X\), gives
\[
 \|x\|_X\leq\|Jx\|_a,
 \qquad
 \|Jx\|_a^a\leq2^a\|x\|_X^a+2^a\|ix\|_X^a
                    =2^{a+1}\|x\|_X^a.
\]
Thus \(J\) is injective and
\(\operatorname{dist}(J)\leq2^{1+1/a}=D_a\).
Its target has \(N\) complex coordinates. The factor two in the
preliminary coordinate bounds is absorbed into the constants stated
in the lemma. This proves the complex case.
\end{proof}

\section{Random compression of finite-dimensional subspaces}\label{sec:compression}
The purpose of this section is to obtain a single map which is
bounded above on the ambient coordinate space and bounded below on
a prescribed \(d\)-dimensional subspace. We first prove an expectation
bound and a concentration estimate for each fixed vector. A covering
argument then yields the upper operator bound. A separate Fourier
argument gives the small-ball estimate needed for the lower bound.
The final covering argument compares the real dimension \(\sigma d\)
of the source subspace with the real dimension \(\sigma m\) of the target.
The range \(a<1\) also makes the defining series absolutely
summable. We can therefore work directly with its characteristic
function, while the covering and perturbation arguments use the
\(a\)- and \(r\)-triangle inequalities. No convexity of the source
or target unit ball is required in these steps.

Fix positive integers \(N,m\), \(0<a<1\), and \(a/2<r<a\),
and write \(\eta=1-r/a\in(0,1/2)\).
Put \(\mathcal U_{\R}=\{1,-1\}\) and
\(\mathcal U_{\mathbb C}=\{1,-1,i,-i\}\).
Let \(e_1,\ldots,e_m\) be the standard coordinate vectors of
\(\K^m\). All random variables in this section are defined on the
product probability space supporting the independent choices below.
The symbols \(\Prb\) and \(\E\) denote its probability and expectation.
The notation \(\|\cdot\|_2\) refers to the Euclidean norm.
Let \(Y_{ij}\) be independent vectors uniform on
\(\{u e_\ell:u\in\mathcal U_{\K},\ 1\leq\ell\leq m\}\subset\K^m\).
Define the \(\K\)-linear map
\begin{equation}\label{eq:operator}
 \mathcal{T}:\ell_a^N\longrightarrow\ell_r^m,\qquad
 \mathcal{T}x=m^{-(1/r-1/a)}\sum_{i=1}^N x_i\sum_{j\geq1}j^{-1/a}Y_{ij}.
\end{equation}
For each outcome, every column series converges absolutely, since
\(\|Y_{ij}\|_2=1\) and \(\sum_jj^{-1/a}<\infty\).
The finitely many column series therefore define \(\mathcal{T}\) simultaneously
for all \(x\), and its entries are measurable.
In particular,
\[
 \|\mathcal{T}x\|_2\leq m^{-(1/r-1/a)}
      \left(\sum_i|x_i|\right)\sum_{j\geq1}j^{-1/a}.
\]
Equivalence of finite-dimensional coordinate quasi-norms gives a
deterministic bound for \(\|\mathcal{T}x\|_r\) at fixed \(N,m,x\).
This is the deterministic-weight construction of \cite{FG11}, up to
a scalar normalization. Related tail estimates appear in \cite{BLV07}.
The estimates below also apply when \(m<N\).
We use this existing construction, not a new distribution of random
operators. The fixed-vector estimates are proved below in the form
needed for the subspace argument. The upper bound is imposed on
the ambient space, whereas the final net is taken only in \(E\).
We prove all estimates used here for the displayed operator, including
the complex case, rather than applying an embedding theorem from
these references to a different range of dimensions. The restriction
\(r<a\) ensures convergence of the moment integral in
Lemma~\ref{lem:mean}. The restriction \(r>a/2\) ensures square
summability of the martingale-difference bounds in
Lemma~\ref{lem:concentration}.

\begin{lemma}\label{lem:mean}
For every unit vector \(x\in\ell_a^N\),
\(\E\|\mathcal{T}x\|_r^r\leq M_{a,r}<\infty\), uniformly in \(N,m,x\).
\end{lemma}
\begin{proof}
Discard the zero coefficients and rearrange the positive terms
\((|x_i|j^{-1/a})_{i,j}\) in nonincreasing order as \((\beta_k)\),
breaking ties deterministically. Absorb each nonzero coefficient's phase into its
random vector and denote the resulting vectors by \(Z_k\).
They remain independent, with \(\E Z_k=0\),
\(\Prb\{Z_{k,\ell}\ne0\}=1/m\), and
\(\E|Z_{k,\ell}|^2=1/m\). The phase absorbed from a nonzero
scalar \(x_i\) is \(x_i/|x_i|\).
The argument uses only these properties of the rotated vectors.
More explicitly, for a fixed \(x\), both the ordering and each
absorbed phase are deterministic. Thus independence is preserved,
although the complex rotations need not have the original four-point
phase distribution. The series can be reordered because
\(\sum_{i,j}|x_i|j^{-1/a}<\infty\).
For \(u>0\), let \(n(u)=\#\{k:\beta_k>u\}\), where \(\#\)
denotes the number of indices in the set. Then
\begin{align*}
 n(u)
 &=\sum_{i=1}^N\#\{j\geq1:|x_i|j^{-1/a}>u\}\\
 &=\sum_{i=1}^N\#\{j\geq1:j<(|x_i|/u)^a\}
 \leq u^{-a}\sum_{i=1}^N|x_i|^a=u^{-a}.
\end{align*}
For \(0<u<\beta_k\), we have \(k\leq n(u)\leq u^{-a}\).
Letting \(u\uparrow\beta_k\) yields \(\beta_k\leq k^{-1/a}\).
Put \(W=\sum_k\beta_k Z_{k,1}\), the first coordinate of the
unnormalized series. The probability that a term with \(\beta_k>t\)
has nonzero first coordinate is at most \(m^{-1}t^{-a}\).
Write \(W_{\leq t}=\sum_{\beta_k\leq t}\beta_kZ_{k,1}\).
Since \(\beta_k\leq k^{-1/a}\) and \(2/a>1\),
\[
 \sum_k\beta_k^2\leq\sum_k k^{-2/a}<\infty.
\]
Independence and centering imply that the partial sums defining
\(W_{\leq t}\) converge in \(L_2\). Their limit agrees with the
absolutely convergent series, and passage to the limit in the
second-moment identity gives
\[
 \E|W_{\leq t}|^2=\frac1m\sum_{\beta_k\leq t}\beta_k^2.
\]
The square sum satisfies
\begin{align*}
 \sum_{\beta_k\leq t}\beta_k^2
 &=\int_0^t2u\,\#\{k:u<\beta_k\leq t\}\,du\\
 &\leq\int_0^t2u\,n(u)\,du
 \leq2\int_0^t u^{1-a}\,du
 =\frac{2}{2-a}t^{2-a}.
\end{align*}
Let \(\mathcal A_t=\{\exists k:\beta_k>t,\ Z_{k,1}\ne0\}\).
On \(\mathcal A_t^c\), we have \(W=W_{\leq t}\). Thus
\begin{align*}
 \Prb\{|W|>t\}
 &\leq\Prb(\mathcal A_t)+\Prb\{|W_{\leq t}|>t\}\\
 &\leq\frac{n(t)}m+t^{-2}\E|W_{\leq t}|^2
 \leq\left(1+\frac2{2-a}\right)\frac1{mt^a}.
\end{align*}
Taking \(C_a=1+2/(2-a)\), this is
\(\Prb\{|W|>t\}\leq\min\{1,C_a/(mt^a)\}\).
Set \(t_0=(C_a/m)^{1/a}\).
Integrating the tail bound and using \(r<a\) gives
\begin{align*}
 \E|W|^r
 &=r\int_0^\infty t^{r-1}\Prb\{|W|>t\}\,dt\\
 &\leq r\int_0^{t_0}t^{r-1}\,dt
   +\frac{rC_a}{m}\int_{t_0}^{\infty}t^{r-a-1}\,dt\\
 &=t_0^r+\frac{rC_a}{m(a-r)}t_0^{r-a}
 =\frac{a}{a-r}(C_a/m)^{r/a}.
\end{align*}
For complex scalars, the same second-moment identity holds because
independence and centering give
\(\E(Z_{k,1}\overline{Z_{l,1}})=0\) for \(k\ne l\).
For \(1\leq\ell\leq m\), put
\(W_\ell=\sum_k\beta_kZ_{k,\ell}\). The same calculation gives
\(\E|W_\ell|^r\leq C_{a,r}m^{-r/a}\). Therefore
\begin{align*}
 \E\|\mathcal{T}x\|_r^r
 &=m^{-r(1/r-1/a)}\sum_{\ell=1}^m\E|W_\ell|^r\\
 &\leq m^{-(1-r/a)}\,m\,C_{a,r}m^{-r/a}
 =C_{a,r}.
\end{align*}
\end{proof}

\begin{lemma}\label{lem:concentration}
There is \(c_{a,r}>0\) such that for a unit \(x\) and every \(t>0\),
\[
 \Prb\{|\|\mathcal{T}x\|_r^r-\E\|\mathcal{T}x\|_r^r|>t\}
 \leq 2\exp(-c_{a,r}m t^{1/\eta}).
\]
\end{lemma}
\begin{proof}
Use the rearrangement and the vectors \(Z_k\) from
Lemma~\ref{lem:mean}. Replacing \(Z_k\) by another possible value, while keeping all other
vectors fixed, changes \(F=\|\mathcal{T}x\|_r^r\) by at most
\(2m^{-\eta}\beta_k^r\). Let \(\mathcal F_k=\sigma(Z_1,\ldots,Z_k)\),
\(M_k=\E(F\mid\mathcal F_k)\), and \(d_k=M_k-M_{k-1}\),
with \(M_0=\E F\). This is the Doob martingale of \(F\).
Put \(\mathcal F_\infty=\sigma(Z_1,Z_2,\ldots)\).
To see why the same bound applies to
\(d_k\), fix \(Z_1,\ldots,Z_{k-1}\) and average \(F\) over the
independent variables with indices greater than \(k\). As a function
of \(Z_k\), this conditional average has oscillation at most
\(2m^{-\eta}\beta_k^r\). Its mean over \(Z_k\) is \(M_{k-1}\),
so its distance from that mean is at most the same oscillation. Hence
\[
 |d_k|\leq 2m^{-\eta}k^{-r/a}.
\]
Indeed, if \(\mathcal{T}'x\) is the output after replacing \(Z_k\) by
another possible value \(Z'_k\), then
\[
 \bigl|\|\mathcal{T}x\|_r^r-\|\mathcal{T}'x\|_r^r\bigr|
 \leq\|(\mathcal{T}-\mathcal{T}')x\|_r^r
 =m^{-\eta}\beta_k^r\|Z_k-Z'_k\|_r^r
 \leq2m^{-\eta}\beta_k^r.
\]
For fixed \(N,m,x\), absolute convergence of the defining series
implies that \(F\) is bounded. Its Doob martingale therefore converges
almost surely and in \(L_2\) to \(F\), since \(F\) is
\(\mathcal F_\infty\)-measurable.

We separate the first \(K\) differences, which will be bounded
deterministically, from the remaining differences, to which an
exponential martingale estimate will be applied. Choosing \(K\)
as a function of the deviation \(t\) produces the power
\(t^{1/\eta}\) in the conclusion.
Write \(\alpha=r/a=1-\eta>1/2\).
For \(h=\eta t m^\eta/4\geq1\), choose
\(K=\lfloor h^{1/\eta}\rfloor\). Since \(1-\alpha=\eta\),
\[
 \sum_{k=1}^K k^{-\alpha}
 \leq1+\int_1^K u^{-\alpha}\,du
 =1+\frac{K^\eta-1}{\eta}
 \leq\frac{K^\eta}{\eta},
 \qquad
 \sum_{k>K}k^{-2\alpha}
 \leq\int_K^\infty u^{-2\alpha}\,du
 =\frac{K^{1-2\alpha}}{2\alpha-1}.
\]
Consequently,
\[
 \sum_{k\leq K}|d_k|
 \leq \frac{2}{\eta}m^{-\eta}K^\eta\leq t/2,
 \qquad
 \sum_{k>K}(2m^{-\eta}k^{-\alpha})^2
 \leq\frac4{2\alpha-1}m^{-2\eta}K^{1-2\alpha}.
\]
Put \(\gamma_k=2m^{-\eta}k^{-\alpha}\) and
\[
 V_K=\frac4{2\alpha-1}m^{-2\eta}K^{1-2\alpha}.
\]
For \(n>K\),
\[
 M_n-M_K=\sum_{k=K+1}^n d_k,
 \qquad
 \sum_{k=K+1}^n\gamma_k^2\leq V_K.
\]
The Azuma--Hoeffding inequality gives, for \(s>0\),
\[
 \Prb\{|M_n-M_K|>s\}
 \leq2\exp\left(-\frac{s^2}{2V_K}\right).
\]
By almost sure convergence, apply this bound with \(s<t/2\),
let \(n\to\infty\), and then let \(s\uparrow t/2\). We obtain
\[
 \Prb\{|F-M_K|>t/2\}
 \leq2\exp\left(-\frac{t^2}{8V_K}\right).
\]
Moreover,
\[
 |M_K-\E F|\leq\sum_{k=1}^K|d_k|\leq t/2,
\]
so
\begin{align*}
 \Prb\{|F-\E F|>t\}
 &\leq\Prb\{|F-M_K|>t/2\}\\
 &\leq2\exp\left(-\frac{t^2}{8V_K}\right)\\
 &=2\exp\left(-\frac{2\alpha-1}{32}
                 t^2m^{2\eta}K^{2\alpha-1}\right).
\end{align*}
Since \(h\geq1\),
\[
 K=\lfloor h^{1/\eta}\rfloor\geq\tfrac12h^{1/\eta}.
\]
Using \(h=\eta tm^\eta/4\) and \(\alpha=1-\eta\), we obtain
\begin{align*}
 t^2m^{2\eta}K^{2\alpha-1}
 &\geq2^{-(2\alpha-1)}t^2m^{2\eta}h^{(2\alpha-1)/\eta}\\
 &=2^{-(2\alpha-1)}
   \left(\frac\eta4\right)^{(2\alpha-1)/\eta}
   t^{2+(2\alpha-1)/\eta}m^{2\eta+2\alpha-1}\\
 &=\kappa_\eta\,m t^{1/\eta},
\end{align*}
where
\[
 \kappa_\eta=2^{-(2\alpha-1)}
       \left(\frac\eta4\right)^{(2\alpha-1)/\eta}>0.
\]
Thus, for \(h\geq1\),
\[
 \Prb\{|F-\E F|>t\}
 \leq2\exp\left(-\frac{2\alpha-1}{32}
                      \kappa_\eta m t^{1/\eta}\right).
\]
For \(h<1\),
\[
 m t^{1/\eta}<\left(\frac4\eta\right)^{1/\eta}.
\]
Choose
\[
 c_{a,r}=\min\left\{
   \frac{2\alpha-1}{32}\kappa_\eta,
   (\log2)\left(\frac\eta4\right)^{1/\eta}\right\}>0.
\]
Then, when \(h<1\),
\[
 \Prb\{|F-\E F|>t\}\leq1
 \leq2\exp(-c_{a,r}m t^{1/\eta}).
\]
The same bound holds when \(h\geq1\) by the preceding estimate,
which completes the proof.
\end{proof}

We shall use the following elementary covering estimate.
For a finite set \(\mathcal N\), \(|\mathcal N|\) or
\(\#\mathcal N\) denotes its cardinality. For a measurable body
\(B\) in a finite-dimensional real space, \(|B|\) denotes Lebesgue
volume and \(\operatorname{int}B\) its interior.

The unit sphere of a \(d\)-dimensional \(a\)-normed space over
\(\K\) admits a \(\delta\)-net, measured in its \(a\)-norm,
of cardinality at most
\(3^{\sigma d/a}\delta^{-\sigma d}\), \(0<\delta\leq1\).
Indeed, let \(B\) be the closed unit ball. For any finite subset
\(\mathcal N\) of the unit sphere whose distinct points have distance
greater than \(\delta\), the translates
\(y+(\delta/2^{1/a})\operatorname{int}B\), \(y\in\mathcal N\),
are disjoint and lie in \((1+\delta^a/2)^{1/a}B\).
Indeed, a point in the intersection of the translates centered at
\(y\) and \(z\) would imply
\(\|y-z\|^a<\delta^a/2+\delta^a/2=\delta^a\), a contradiction.
For a point \(y+u\) in one translate,
\(\|y+u\|^a\leq1+\|u\|^a<1+\delta^a/2\), which proves
containment. This argument uses the \(a\)-triangle inequality and
homogeneity, and does not require convexity of \(B\).
Since \(sB\subset\operatorname{int}B\) for \(0<s<1\),
\(|\operatorname{int}B|\geq s^{\sigma d}|B|\) for every \(s<1\).
Letting \(s\uparrow1\) gives \(|\operatorname{int}B|=|B|\).
The unit ball is compact and has positive finite volume, since its
\(a\)-norm is equivalent to a Euclidean norm in finite dimension.
Lebesgue volume in real dimension
\(\sigma d\) therefore gives
\[
 |\mathcal N|\left(\frac\delta{2^{1/a}}\right)^{\sigma d}|B|
 \leq(1+\delta^a/2)^{\sigma d/a}|B|.
\]
Consequently,
\[
 |\mathcal N|\leq(1+2/\delta^a)^{\sigma d/a}
 \leq3^{\sigma d/a}\delta^{-\sigma d}.
\]
Starting from one point, successively add a sphere point at distance
greater than \(\delta\) from all previously selected points whenever
possible. The cardinality bound forces this procedure to terminate.
The resulting set is a \(\delta\)-net by construction.

All random operator extrema used below are measurable. Indeed, for any
fixed nonzero subspace \(E\subset\ell_a^N\), choose a countable dense
subset \(\mathcal D_E\) of its unit sphere. Continuity gives
\[
 \sup_{\|x\|_a=1,\,x\in E}\|\mathcal{T}x\|_r
   =\sup_{x\in\mathcal D_E}\|\mathcal{T}x\|_r,
 \qquad
 \inf_{\|x\|_a=1,\,x\in E}\|\mathcal{T}x\|_r
   =\inf_{x\in\mathcal D_E}\|\mathcal{T}x\|_r.
\]
The right-hand sides are countable extrema of measurable functions.

\begin{lemma}\label{lem:upper}
For each \(C\geq1\) there is \(L=L(a,r,C)\geq1\) such that
\[
 N/m\leq C
 \quad\Longrightarrow\quad
 \Prb\{\|\mathcal{T}:\ell_a^N\to\ell_r^m\|\leq L\}\geq3/4.
\]
\end{lemma}
\begin{proof}
Put \(\delta=2^{-1/r}\), and choose a \(\delta\)-net
\(\mathcal N\) of the source unit sphere. The covering estimate gives
\[
 |\mathcal N|\leq H^{\sigma N},\qquad H=3^{1/a}2^{1/r}.
\]
For each unit vector \(x\), choose \(y\in\mathcal N\) with
\(\|x-y\|_a\leq\delta\). Then
\[
 \|\mathcal{T}x\|_r^r\leq\|\mathcal{T}y\|_r^r+\|\mathcal{T}(x-y)\|_r^r
 \leq\max_{z\in\mathcal N}\|\mathcal{T}z\|_r^r+\delta^r\|\mathcal{T}\|^r.
\]
Since \(\|\mathcal{T}\|<\infty\) and \(\delta^r=1/2\), taking the supremum yields
\[
 \tfrac12\|\mathcal{T}\|^r\leq\max_{z\in\mathcal N}\|\mathcal{T}z\|_r^r.
\]
The net is deterministic, so the fixed-vector estimates apply to each
of its points. For every \(z\in\mathcal N\),
\(\E\|\mathcal{T}z\|_r^r\leq M_{a,r}\). Hence the event
\(\|\mathcal{T}z\|_r^r>M_{a,r}+t\) implies a deviation of more than \(t\)
above its own expectation. A union bound requires no independence
between the events at different net points.
By Lemmas~\ref{lem:mean} and \ref{lem:concentration}, for \(t>0\),
\begin{align*}
 \Prb\{\|\mathcal{T}\|^r>2(M_{a,r}+t)\}
 &\leq\sum_{z\in\mathcal N}
       \Prb\{\|\mathcal{T}z\|_r^r>M_{a,r}+t\}\\
 &\leq2|\mathcal N|e^{-c_{a,r}mt^{1/\eta}}\\
 &\leq2\exp\bigl(\sigma N\log H-c_{a,r}mt^{1/\eta}\bigr).
\end{align*}
Choose \(t>0\) such that
\(c_{a,r}t^{1/\eta}\geq2C\log H+\log8\). Using \(N\leq Cm\)
and \(\sigma\leq2\), we obtain
\[
 \sigma N\log H-c_{a,r}mt^{1/\eta}
 \leq2Cm\log H-m(2C\log H+\log8)=-m\log8.
\]
With \(L=\max\{1,[2(M_{a,r}+t)]^{1/r}\}\), it follows that
\[
 \Prb\{\|\mathcal{T}\|>L\}\leq2\cdot8^{-m}\leq\tfrac14.
\]
\end{proof}

A small-ball estimate bounds the probability that the image of a
fixed unit vector lies in a small ball about zero. The next lemma
gives such a bound with the power \(\sigma m\), the real dimension
of the target. This precise power is needed: when the bound is
summed over a net in a \(d\)-dimensional subspace, the net contributes
\(t^{-\sigma d}\), leaving \(t^{\sigma(m-d)}\). The hypothesis
\(N/m\leq C\) controls the constants even when \(m<N\). It does
not assert injectivity on the whole ambient space.

\begin{lemma}\label{lem:density}
If \(C\geq1\) and \(N/m\leq C\), there is \(B=B(a,r,C)\geq1\) such that
for every unit vector \(x\in\ell_a^N\) and every \(t>0\),
\[
 \Prb\{\|\mathcal{T}x\|_r\leq t\}\leq (Bt)^{\sigma m}.
\]
\end{lemma}
\begin{proof}
Fix a unit vector \(x\). We will bound the density of the rescaled
random vector \(Z=m^{1/r}\mathcal{T}x\), then multiply that bound by the
volume of a ball. The rescaling is chosen so that the factor
\(m^{-1/r}\) in the volume estimate cancels exactly. Although the
coordinates of \(Z\) are not asserted to be independent, the
independent summands give a product bound for its characteristic
function.
For \(z\in\K\), let \(z_1,\ldots,z_\sigma\) be its real
components and put
\[
 \Delta_{\K}(z)^2=
   \sum_{h=1}^{\sigma}\operatorname{dist}(z_h,\pi\mathbb Z)^2.
\]
In this formula,
\(\operatorname{dist}(u,\pi\mathbb Z)=\inf_{n\in\mathbb Z}|u-\pi n|\)
is scalar distance to the lattice \(\pi\mathbb Z\). It is distinct
from the distortion of a linear map.
Set
\[
 Z_J=m^{1/a}\sum_{i=1}^N\sum_{j=1}^Jx_i j^{-1/a}Y_{ij},
 \qquad Z=m^{1/r}\mathcal{T}x=\lim_{J\to\infty}Z_J.
\]
The limit holds for every outcome by absolute convergence. Use the real
pairing \([\xi,u]=\operatorname{Re}\sum_\ell\overline{\xi_\ell}u_\ell\).
For each \((i,j)\), put
\[
 \phi_{ij}(\xi)
 =\E e^{i[\xi,m^{1/a}x_i j^{-1/a}Y_{ij}]}
 =\frac1{\sigma m}\sum_{\ell=1}^m\sum_{h=1}^{\sigma}
 \cos\bigl((m^{1/a}x_i\overline{\xi_\ell})_h j^{-1/a}\bigr).
\]
For real scalars, averaging the signs gives the cosine in this
formula. For complex scalars, averaging the four phases gives
\((\cos(\operatorname{Re}z)+\cos(\operatorname{Im}z))/2\) at
each coordinate, which explains the factor \(1/(\sigma m)\).
If \(u=\operatorname{dist}(v,\pi\mathbb Z)\in[0,\pi/2]\),
then \(|\cos v|=\cos u\) and
\(1-\cos u=2\sin^2(u/2)\geq2u^2/\pi^2\).
Thus an absolute constant \(c>0\) works uniformly in \(v\).
The bounds \(|\cos v|\leq1-c\operatorname{dist}(v,\pi\mathbb Z)^2\)
and \(1-u\leq e^{-u}\) give
\begin{align*}
 |\phi_{ij}(\xi)|
 &\leq1-\frac{c}{\sigma m}\sum_{\ell=1}^m
   \Delta_{\K}(m^{1/a}x_i\overline{\xi_\ell}j^{-1/a})^2\\
 &\leq\exp\left(-\frac{c}{\sigma m}\sum_{\ell=1}^m
   \Delta_{\K}(m^{1/a}x_i\overline{\xi_\ell}j^{-1/a})^2\right).
\end{align*}
By independence and \(\sigma\leq2\), after decreasing \(c\),
\begin{align*}
 |\E e^{i[\xi,Z_J]}|
 &=\prod_{i=1}^N\prod_{j=1}^J|\phi_{ij}(\xi)|\\
 &\leq\prod_{\ell=1}^m
 \exp\left(-\frac c m\sum_{i=1}^N\sum_{j=1}^J
 \Delta_{\K}(m^{1/a}x_i\overline{\xi_\ell}j^{-1/a})^2\right).
\end{align*}
For fixed \(A\in\K\), the summands
\(\Delta_{\K}(Aj^{-1/a})^2\) equal \(|A|^2j^{-2/a}\) for all
sufficiently large \(j\). Since \(2/a>1\), their series converges.
Bounded convergence on the left and convergence of these series on
the right therefore yield
\begin{equation}\label{eq:characteristic-product}
 |\E e^{i[\xi,Z]}|\leq\prod_{\ell=1}^m H_x(\xi_\ell),
 \qquad
 H_x(z)=\exp\left(-\frac c m\sum_{i=1}^N\sum_{j\geq1}
 \Delta_{\K}(m^{1/a}x_i\overline z j^{-1/a})^2\right).
\end{equation}
For \(A\in\K\), consider the terms with
\(j\geq\max\{1,\lceil|A|^a\rceil\}\). Each real component then
has absolute value at most one, and
\(\Delta_{\K}(Aj^{-1/a})^2=|A|^2j^{-2/a}\).
For \(|A|\leq1\), the term \(j=1\) alone is \(|A|^2\).
For \(|A|>1\), put \(J=\lceil|A|^a\rceil\leq2|A|^a\).
Then
\[
 |A|^2\sum_{j\geq J}j^{-2/a}
 \geq\frac{a}{2-a}|A|^2J^{1-2/a}
 \geq c_a|A|^a.
\]
These two cases give
\[
 \sum_{j\geq1}\Delta_{\K}(Aj^{-1/a})^2
 \geq c_a\min\{|A|^2,|A|^a\}.
\]
Put \(v_i=|x_i|^a|z|^a\), so that \(\sum_i v_i=|z|^a\).
It follows that
\[
 -\log H_x(z)\geq c_a\sum_i
       \min\{m^{2/a-1}v_i^{2/a},v_i\}.
\]
To check the powers of \(m\), substitute
\(A=m^{1/a}x_i\overline z\) in the preceding estimate and divide
by \(m\). Its two terms become
\(m^{-1}|A|^2=m^{2/a-1}v_i^{2/a}\) and
\(m^{-1}|A|^a=v_i\), respectively.
For indices satisfying \(mv_i\geq1\), the minimum equals \(v_i\).
If the sum of \(v_i\) over these indices is at least half the total,
the sum on the right-hand side is at least
\(|z|^a/2\). Otherwise, power-mean comparison on the remaining
indices yields
\[
 m^{2/a-1}\sum_{mv_i<1}v_i^{2/a}
 \geq m^{2/a-1}N^{1-2/a}(|z|^a/2)^{2/a}
 \geq2^{-2/a}C^{1-2/a}|z|^2.
\]
For completeness, if \(J_0=\{i:mv_i<1\}\) has \(n_0\leq N\)
elements and its sum is positive, H\"older's inequality gives
\[
 \sum_{i\in J_0}v_i
 \leq n_0^{1-a/2}
       \left(\sum_{i\in J_0}v_i^{2/a}\right)^{a/2}.
\]
Rearranging yields the first lower bound above, because
\(n_0^{1-2/a}\geq N^{1-2/a}\). Moreover,
\[
 m^{2/a-1}N^{1-2/a}=(N/m)^{1-2/a}\geq C^{1-2/a}.
\]
The last inequality uses \(1-2/a<0\). Hence, uniformly in \(x,N,m\),
\begin{equation}\label{eq:uniform-fourier-envelope}
 H_x(z)\leq
 \exp\left(-c_a'\min\{|z|^a,C^{1-2/a}|z|^2\}\right)
 =:H_{a,C}(z).
\end{equation}
For \(|z|\geq C^{1/a}\),
\[
 C^{1-2/a}|z|^2\geq|z|^a,
 \qquad H_{a,C}(z)=e^{-c_a'|z|^a}.
\]
Thus \(I_{a,C,\sigma}:=\int_{\K}H_{a,C}(z)\,dz<\infty\), and
\begin{align*}
 \int_{\K^m}|\E e^{i[\xi,Z]}|\,d\xi
 &\leq\int_{\K^m}\prod_{\ell=1}^mH_{a,C}(\xi_\ell)\,d\xi\\
 &=I_{a,C,\sigma}^m<\infty.
\end{align*}
Fourier inversion therefore gives a bounded continuous joint
probability density with respect to Lebesgue measure on
\(\K^m\cong\R^{\sigma m}\), namely
\[
 f_Z(u)=(2\pi)^{-\sigma m}
       \int_{\K^m}e^{-i[\xi,u]}\E e^{i[\xi,Z]}\,d\xi.
\]
Taking absolute values yields
\[
 \|f_Z\|_\infty\leq
 \bigl((2\pi)^{-\sigma}I_{a,C,\sigma}\bigr)^m
 \leq D_{a,C}^{\sigma m},
\]
where \(D_{a,C}\) is chosen for both \(\sigma=1,2\).
Let \(B_r^m(\K)\) be the unit ball, and put

\(c_{r,\sigma}=\int_{\R^\sigma}e^{-|z|^r}\,dz\) and
\(\lambda=\sigma m/r\). Then
\begin{align*}
 e^{-\lambda}|B_r^m(\K)|
 &\leq\int_{\R^{\sigma m}}
       e^{-\lambda\sum_{\ell=1}^m|u_\ell|^r}\,du\\
 &=\prod_{\ell=1}^m\int_{\R^\sigma}e^{-\lambda|z|^r}\,dz
 =\bigl(c_{r,\sigma}\lambda^{-\sigma/r}\bigr)^m.
\end{align*}
Thus
\[
 |B_r^m(\K)|\leq
 \left[e^{1/r}c_{r,\sigma}^{1/\sigma}
       (r/\sigma)^{1/r}m^{-1/r}\right]^{\sigma m}
 \leq(C_rm^{-1/r})^{\sigma m},
\]
where
\[
 C_r=\max_{\sigma\in\{1,2\}}
 e^{1/r}c_{r,\sigma}^{1/\sigma}(r/\sigma)^{1/r}.
\]
Consequently,
\[
 \Prb\{\|\mathcal{T}x\|_r\leq t\}
 =\Prb\{Z\in tm^{1/r}B_r^m(\K)\}
 \leq\|f_Z\|_\infty(tm^{1/r})^{\sigma m}|B_r^m(\K)|
 \leq(D_{a,C}C_rt)^{\sigma m}.
\]
The conclusion follows with \(B=\max\{1,D_{a,C}C_r\}\).
\end{proof}

\begin{proposition}\label{prop:compression}
Let \(0<a<1\), \(a/2<r<a\), \(C_0\geq1\), and \(\mu>1\).
Every \(d\)-dimensional \(E\subset\ell_a^N\) with \(N\leq C_0d\)
embeds into \(\ell_r^{\lfloor\mu d\rfloor}\) with distortion at most
\(K(a,r,C_0,\mu)\).
\end{proposition}
\begin{proof}
The case \(d=0\) is immediate. Assume \(d\geq1\).
Put \(m=\lfloor\mu d\rfloor\). Then \(m\geq d\) and \(N/m\leq C_0\).
Consider the random operator \eqref{eq:operator}, with \(L\) and \(B\)
as in Lemmas~\ref{lem:upper} and \ref{lem:density}.
For \(0<t\leq1\), choose a
\(\delta=t/(2^{1/r}L)\)-net of the sphere of \(E\),
chosen deterministically before sampling \(\mathcal{T}\).
Write \(\mathcal N_E\) for this net and define
\[
 \mathcal G_t=\left\{\|\mathcal{T}\|\leq L,\quad
 \inf_{\substack{x\in E\\\|x\|_a=1}}\|\mathcal{T}x\|_r
 <t/2^{1/r}\right\}.
\]
If \(\|\mathcal{T}\|\leq L\) and a unit vector \(x\in E\) satisfies
\(\|\mathcal{T}x\|_r<t/2^{1/r}\), a net point \(y\) within distance \(\delta\) satisfies
\(\|\mathcal{T}y\|_r<t\), because
\[
 \|\mathcal{T}y\|_r^r\leq\|\mathcal{T}x\|_r^r+\|\mathcal{T}(y-x)\|_r^r
       <t^r/2+L^r\delta^r=t^r.
\]
The union bound and Lemma~\ref{lem:density} give
\begin{align*}
 \Prb(\mathcal G_t)
 &\leq\sum_{y\in\mathcal N_E}\Prb\{\|\mathcal{T}y\|_r<t\}
 \leq|\mathcal N_E|(Bt)^{\sigma m}\\
 &\leq3^{\sigma d/a}
       \left(\frac{2^{1/r}L}{t}\right)^{\sigma d}(Bt)^{\sigma m}
 =\left[(3^{1/a}2^{1/r}L)^d B^m t^{m-d}\right]^\sigma.
\end{align*}
The estimate combines \(t^{-\sigma d}\) from covering the source
sphere with \(t^{\sigma m}\) from the small-ball probability.
Consequently an excess of target dimension over source dimension
allows the failure probability to be made small by choosing \(t\)
small. All constants in this choice will be independent of \(d\).
The exponent of \(t\) is therefore \(\sigma(m-d)\) over either field.
For \(d\geq2/(\mu-1)\) we have
\(m-d\geq\varepsilon d\), \(\varepsilon=(\mu-1)/2\), and \(m\leq\mu d\).
Indeed,
\[
 m-d\geq(\mu-1)d-1\geq\tfrac12(\mu-1)d.
\]
Since \(B\geq1\) and \(0<t\leq1\), these inequalities imply
\(B^m\leq B^{\mu d}\) and \(t^{m-d}\leq t^{\varepsilon d}\).
Put \(H_0=3^{1/a}2^{1/r}L\) and choose
\[
 t=\min\{1,(8H_0B^\mu)^{-1/\varepsilon}\}.
\]
Then
\[
 \Prb(\mathcal G_t)
 \leq\bigl[H_0B^\mu t^\varepsilon\bigr]^{\sigma d}
 \leq8^{-\sigma d}\leq1/8.
\]
Consequently,
\[
 \Prb\left\{\|\mathcal{T}\|\leq L,\
 \inf_{x\in E,\,\|x\|_a=1}\|\mathcal{T}x\|_r\geq t/2^{1/r}\right\}
 \geq1-\frac14-\frac18=\frac58>0.
\]
Fix an operator in this event. By homogeneity,
\[
 \frac{t}{2^{1/r}}\|x\|_a\leq\|\mathcal{T}x\|_r\leq L\|x\|_a
 \quad(x\in E),
 \qquad \operatorname{dist}(\mathcal{T}|_E)\leq\frac{2^{1/r}L}{t}.
\]

The remaining positive dimensions satisfy \(d<2/(\mu-1)\).
For these dimensions the floor in \(m=\lfloor\mu d\rfloor\)
can leave no excess coordinates at all. A deterministic coordinate
map suffices, because \(d\) is now bounded in terms of \(\mu\).
The following maximal-determinant argument supplies the needed
coordinate functionals without invoking a convexity assumption.
Fix reference coordinates on \(E\). Compactness of its unit sphere
implies that the absolute determinant of a \(d\)-tuple of unit vectors
attains its maximum. This maximum is positive since \(E\) has a basis.
Let \((e_i)\) be a maximizing tuple. Its determinant is nonzero,
so it is a basis. Define the coefficient functionals \(e_i^*\) by
the unique expansion \(x=\sum_{i=1}^d e_i^*(x)e_i\). These are
\(\K\)-linear functionals. Replacing column \(i\) by a unit vector \(x\) multiplies
the determinant by \(e_i^*(x)\). Maximality therefore gives
\(|e_i^*(x)|\leq1\), over either field. By homogeneity,
\(|e_i^*(x)|\leq\|x\|_E\) for every \(x\).
Define \(Vx=(e_i^*(x))_{i=1}^d\). Since \(x=\sum_i e_i^*(x)e_i\),
\[
 \|x\|_E^a\leq\sum_{i=1}^d|e_i^*(x)|^a\|e_i\|_E^a
 =\|Vx\|_a^a,
 \qquad
 \|Vx\|_r^r=\sum_{i=1}^d|e_i^*(x)|^r\leq d\|x\|_E^r.
\]
The coordinate inequality \(\|u\|_a\leq\|u\|_r\) for \(r<a\)
also holds below exponent one: after normalizing \(\|u\|_r=1\),
each \(|u_i|\leq1\), and therefore
\(\sum_i|u_i|^a\leq\sum_i|u_i|^r=1\).
Together with \(r<a\), these inequalities give
\[
 \|x\|_E\leq\|Vx\|_a\leq\|Vx\|_r\leq d^{1/r}\|x\|_E.
\]
The zero-padded map \(\widetilde Vx=(Vx,0)\in\ell_r^m\) satisfies
\[
 \operatorname{dist}(\widetilde V)\leq d^{1/r}
 <\left(\frac2{\mu-1}\right)^{1/r}.
\]
Taking the maximum of the large-dimension bound
\(2^{1/r}L/t\) and the small-dimension bound
\((2/(\mu-1))^{1/r}\) defines \(K(a,r,C_0,\mu)\) and completes
the proof.
\end{proof}

The choice of the small-ball threshold also gives a bound as the
target dimension approaches the source dimension.

\begin{remark}\label{rem:excess}
For \(1<\mu\leq2\), Proposition~\ref{prop:compression} admits the
quantitative choice
\begin{equation}\label{eq:compression-quantitative}
 K(a,r,C_0,\mu)
 \leq \exp\!\left(\frac{C(a,r,C_0)}{\mu-1}\right).
\end{equation}
Indeed, increase \(L\) to be at least one. In the large-dimension
case of the proof, \(H_0\), \(B\), and \(L\) depend only on
\(a,r,C_0\), and
\[
 \log\frac{2^{1/r}L}{t}
 =\log(2^{1/r}L)
   +\frac{2}{\mu-1}\log(8H_0B^\mu)
 \leq\frac{C(a,r,C_0)}{\mu-1}.
\]
For the remaining dimensions, put \(u=\mu-1\in(0,1]\).
The deterministic bound satisfies
\[
 \log\bigl((2/u)^{1/r}\bigr)
 =r^{-1}\log(2/u)\leq\frac{\log2+1}{ru},
\]
which proves \eqref{eq:compression-quantitative} in this case as
well. In Section~\ref{sec:mainproof}, we combine this estimate
with the dependence on \(DT\) in Lemma~\ref{lem:restore} to obtain
\eqref{eq:main-quantitative}.
No optimality of \eqref{eq:compression-quantitative} is asserted.
\end{remark}

The restriction \(r>a/2\) in the random estimate can be met by
first choosing an intermediate exponent. This gives the following
consequence for any smaller target exponent.

\begin{corollary}\label{cor:lower-exponent}
Let \(0<r<p<1\) and \(\mu>1\). Every \(d\)-dimensional
subspace of \(L_p\), \(d\geq1\), admits a \(\K\)-linear
embedding into \(\ell_r^{\lfloor\mu d\rfloor}\) with distortion
at most \(K(p,r,\mu)\), independently of \(d\).
\end{corollary}
\begin{proof}
Choose \(b=(r+\min\{p,2r\})/2\), so that
\(r<b<p\) and \(b/2<r\). Lemma~\ref{lem:preliminary}
embeds the source into a subspace \(E\subset\ell_b^N\) with
\(N\leq C(b,p)d\) and distortion at most \(D_b\).
Proposition~\ref{prop:compression}, with source exponent \(b\)
and target exponent \(r\), embeds \(E\) into the required
coordinate space. The product of the two distortions depends
only on \(p,r,\mu\), since \(b\) is fixed by \(p,r\).
Both maps are \(\K\)-linear.
\end{proof}

\section{Change of density and exponent restoration}\label{sec:restoration}
We now keep the number of coordinates fixed and change their weights.
The type assumption will first give a maximum estimate for the
image of the original unit ball, which is the hypothesis of the
density theorem in
Section~\ref{sec:preliminaries}. Its weak \(L_q\) conclusion implies
a strong \(L_p\) estimate when \(p<q\). At \(p=q\), the tail
integral is logarithmic, so a second change of weights is used to
bound the largest coordinate. All multipliers are chosen for the
whole space, ensuring that the resulting embeddings remain linear.

The following elementary lemma isolates the finite-atomic step.
It applies to arbitrary positive atom probabilities. In particular,
its constants do not depend on their minimum. Throughout its
statement, weighted strong and weak functionals are taken with
respect to the indicated probability vector.

\begin{lemma}\label{lem:regularization}
Let \(0<r<q<\infty\), let \(k\geq1\) be an integer, and let
\(w_1,\ldots,w_k>0\) satisfy
\(\sum_iw_i=1\). Put \(v_i=(w_i+1/k)/2\), and define the
\(\K\)-linear diagonal map
\[
 (Vf)(i)=(w_i/v_i)^{1/r}f(i),\qquad f\in\K^k.
\]
Then \(\sum_i v_i=1\), \(v_i\geq1/(2k)\), and
\begin{align}
 \|Vf\|_{r;v}&=\|f\|_{r;w},\label{eq:regularization-isometry}\\
 \|Vf\|_{q,\infty;v}&\leq C(r,q)\|f\|_{q,\infty;w},
       \label{eq:regularization-weak}\\
 \|Vf\|_{q;v}&\leq(1+\log(2k))^{1/q}\|Vf\|_{q,\infty;v}.
       \label{eq:regularization-strong}
\end{align}
The weights and the map \(V\) are fixed for all \(f\), and the
number of coordinates remains \(k\).
\end{lemma}
\begin{proof}
Put \(v_i=(w_i+1/k)/2\), and define
\(Vf(i)=(w_i/v_i)^{1/r}f(i)\).
The identity
\[
 \sum_i v_i|Vf(i)|^r=\sum_i w_i|f(i)|^r
\]
shows that \(V\) preserves the \(L_r\) functional. The new atom
probabilities satisfy \(v_i\geq1/(2k)\).
They also sum to one, because \(\sum_iw_i=1\). Moreover
\(v_i/w_i=(1+1/(kw_i))/2\geq1/2\), although this ratio need
not have a uniform upper bound. The dyadic decomposition below
handles this possible lack of an upper bound.
For these weighted coordinate spaces the notation is explicitly
\[
 \|f\|_{q;v}=\left(\sum_{i=1}^k v_i|f(i)|^q\right)^{1/q},
 \qquad
 \|f\|_{q,\infty;v}
 =\sup_{t>0}t\left(\sum_{\{i:|f(i)|>t\}}v_i\right)^{1/q}.
\]
We use analogous notation for the probabilities \((w_i)\) and
for other strong exponents, such as \(r\).
The map \(V\) preserves the weak \(L_q\) estimate up to a factor
depending only on \(r,q\). To prove this, set \(a_i=v_i/w_i\geq1/2\) and partition the
indices into \(I_j=\{2^j\leq a_i<2^{j+1}\}\), \(j=-1,0,1,\ldots\).
For \(i\in I_j\), the condition \(|Vf(i)|>t\) implies
\(|f(i)|>t2^{j/r}\), while \(v_i\leq2^{j+1}w_i\).
Thus, if \(\|f\|_{q,\infty;w}\leq A\), summing first over each
\(I_j\) and then dropping that restriction gives
\begin{align*}
 \sum_{\{|Vf(i)|>t\}}v_i
 &\leq\sum_{j\geq-1}2^{j+1}
       \sum_{\{|f(i)|>t2^{j/r}\}}w_i\\
 &\leq 2 A^q t^{-q}\sum_{j\geq-1}2^{j(1-q/r)}.
\end{align*}
Since \(q/r>1\),
\[
 \sum_{j\geq-1}2^{j(1-q/r)}
 =\frac{2^{q/r-1}}{1-2^{1-q/r}}<\infty.
\]
Consequently,
\[
 \|Vf\|_{q,\infty;v}
 \leq\left(\frac{2^{q/r}}{1-2^{1-q/r}}\right)^{1/q}A
 =:C(r,q)A.
\]
Put \(H=\|Vf\|_{q,\infty;v}\). If \(H=0\), then \(Vf=0\).
Assume \(H>0\). For each \(i\) and \(0<s<|Vf(i)|\),
\[
 H\geq s\left(\sum_{\{|Vf(j)|>s\}}v_j\right)^{1/q}
 \geq s v_i^{1/q}\geq s(2k)^{-1/q}.
\]
Letting \(s\uparrow|Vf(i)|\) gives
\[
 \max_i|Vf(i)|\leq(2k)^{1/q}H.
\]
Set \(M=(2k)^{1/q}H\). The tail distribution function is at most one
for \(t<H\), at most \((H/t)^q\) for \(H\leq t\leq M\), and
zero for \(t>M\). These three ranges explain the two integrals
and their upper cutoff. Layer-cake integration gives
\[
 \|Vf\|_{q;v}^q
 \leq q\int_0^H t^{q-1}\,dt
       +qH^q\int_H^M\frac{dt}{t}
 =H^q(1+\log(2k)).
\]
The first identity proves \eqref{eq:regularization-isometry},
the dyadic estimate proves \eqref{eq:regularization-weak}, and the
last integral proves \eqref{eq:regularization-strong}.
\end{proof}

We now restore the exponent by applying the density theorem to
the image of the original unit ball. The regularization in
Lemma~\ref{lem:regularization} is needed only at the type exponent.

\begin{lemma}\label{lem:restore}
Let \(0<r<q\leq2\), \(r<1\), and let \(X\) be a nonzero
finite-dimensional quasi-normed space over \(\K\) with
\(T_q(X)\leq T\). Suppose that a \(\K\)-linear map
\(A:X\to\ell_r^m\) satisfies
\[
 \|x\|_X\leq\|Ax\|_r\leq D\|x\|_X
 \qquad(x\in X).
\]
The following assertions hold with constants independent of
\(m,D,T\), apart from the displayed factors.
\begin{enumerate}
\item If \(r<p<q\), there is a \(\K\)-linear map
\(R:A(X)\to\ell_p^m\) such that
\[
 \|x\|_X\leq\|RAx\|_p
       \leq C(r,p,q)DT\|x\|_X.
\]
\item There is a \(\K\)-linear map \(R:A(X)\to\ell_q^m\)
such that
\[
 \|x\|_X\leq\|RAx\|_q
       \leq C(r,q)DT(1+\log(2m))^{1/q}\|x\|_X.
\]
\end{enumerate}
\end{lemma}
\begin{proof}
Let \(\nu\) be the uniform probability measure on
\(\{1,\ldots,m\}\), and define
\(\widehat A x=m^{1/r}Ax\). Then
\[
 \|\widehat A x\|_{r;\nu}^r
 =m^{-1}\sum_{j=1}^m|m^{1/r}(Ax)_j|^r=\|Ax\|_r^r.
\]
Thus \(\widehat A\) has the same upper and lower bounds as
\(A\), now with a probability measure as ambient measure.

The scalar Khintchine inequality, over either field, gives
\[
 \E\left|\sum_i\varepsilon_i u_i\right|^r
 \geq c_r\left(\sum_i|u_i|^2\right)^{r/2}
 \geq c_r\max_i|u_i|^r
 \qquad(u_i\in\K),
\]
where \(c_r>0\) depends only on \(r\).

We apply the scalar estimate at each atom to
\(u_i=c_i\widehat A x_i(s)\), where \(x_i\in B_X\)
and \(c_i\in[0,1]\). Tonelli's theorem gives
\begin{align*}
 \left\|\max_i|c_i\widehat A x_i|\right\|_{r;\nu}^r
 &\leq c_r^{-1}\E
       \left\|\sum_i\varepsilon_i c_i\widehat A x_i\right\|_{r;\nu}^r\\
 &\leq c_r^{-1}D^r\E
       \left\|\sum_i\varepsilon_i c_i x_i\right\|_X^r.
\end{align*}
Since \(r<2\), comparison of moments on the probability space of
signs and the type assumption yield
\begin{align*}
 \left(\E\left\|\sum_i\varepsilon_i c_i x_i\right\|_X^r\right)^{1/r}
 &\leq\left(\E\left\|\sum_i\varepsilon_i c_i x_i\right\|_X^2\right)^{1/2}\\
 &\leq T\left(\sum_i|c_i|^q\|x_i\|_X^q\right)^{1/q}
 \leq T\left(\sum_i|c_i|^q\right)^{1/q}.
\end{align*}
Consequently,
\[
 \left\|\max_i|c_i\widehat A x_i|\right\|_{r;\nu}
 \leq c_r^{-1/r}DT\left(\sum_i|c_i|^q\right)^{1/q}.
\]
Apply the density theorem from item~\ref{input:density} to the
real nonnegative functions
\[
 S_0=\{|\widehat A x|:x\in B_X\}\subset L_r(\nu;\R).
\]
All its hypotheses hold, with source exponent \(r\), target
exponent \(q\), and maximum constant \(c_r^{-1/r}DT\).
It gives one probability density \(g\) such that, with
\(\rho=g\,d\nu\), the fixed multiplier
\(Uh=g^{-1/r}h\), defined as zero on \(\{g=0\}\), satisfies
\begin{equation}\label{eq:restoration-weak}
 \|U\widehat A x\|_{q,\infty;\rho}
 \leq C_1(r,q)DT\|x\|_X,
 \qquad
 \|U\widehat A x\|_{r;\rho}=\|Ax\|_r.
\end{equation}
Here the weak estimate first holds for \(x\in B_X\) and then
for all \(x\) by homogeneity. Although taking absolute values
is not linear, the resulting map \(U\widehat A\) is
\(\K\)-linear, since its multiplier is fixed and real
nonnegative. The support assertion of the density theorem implies
that every \(\widehat A x\) vanishes at each atom where
\(g=0\). This also proves the isometry identity in
\eqref{eq:restoration-weak} over either scalar field.
Apply a common permutation to the coordinates of \(A\),
\(\widehat A\), and \(g\), so that the positive values of \(g\)
occur at \(1,\ldots,k\). We retain the same notation after this
permutation, which preserves all the preceding estimates.
In the original coordinate order, the restoration map is obtained
by first applying this permutation and then the formula for \(R\)
given below in the relevant case.
Set \(w_i=g(i)/m\) for \(1\leq i\leq k\). These positive
weights sum to one. Every vector in \(A(X)\) vanishes on the
remaining \(m-k\) coordinates. We identify \(U\widehat A x\)
with its restriction to these \(k\) atoms when applying
Lemma~\ref{lem:regularization}.

Suppose first that \(r<p<q\), and put
\(h=U\widehat A x\) and \(F=C_1(r,q)DT\|x\|_X\).
For \(x\ne0\), \eqref{eq:restoration-weak} gives
\(\rho\{|h|>t\}\leq\min\{1,(F/t)^q\}\).
Layer-cake integration therefore yields
\begin{align*}
 \|h\|_{p;\rho}^p
 &=p\int_0^\infty t^{p-1}\rho\{|h|>t\}\,dt\\
 &\leq p\int_0^F t^{p-1}\,dt
       +pF^q\int_F^\infty t^{p-q-1}\,dt
 =\frac q{q-p}F^p.
\end{align*}
The tail integral converges because \(p<q\). Since \(\rho\)
is a probability and \(r<p\), H\"older's inequality gives
\[
 \|x\|_X\leq\|Ax\|_r=\|h\|_{r;\rho}
 \leq\|h\|_{p;\rho}
 \leq\left(\frac q{q-p}\right)^{1/p}C_1(r,q)DT\|x\|_X.
\]
For \(y\in A(X)\), define
\[
 Ry=\bigl(w_i^{1/p}[U(m^{1/r}y)](i)\bigr)_{i=1}^k
          \oplus0\in\ell_p^m.
\]
The factor \(m^{1/r}\) identifies the original coordinate space
with \(L_r(\nu)\). The density multiplier and the remaining
coordinate weights are fixed, so \(R\) is \(\K\)-linear.
Appending \(m-k\) zero coordinates leaves the functional
unchanged. Since \(m^{1/r}Ax=\widehat A x\), we have
\(\|RAx\|_p=\|h\|_{p;\rho}\). The preceding inequalities
prove (i). They also show that \(R\) is injective on \(A(X)\).
The same inequalities hold trivially for \(x=0\).

At the endpoint \(p=q\), the preceding tail integral would
contain \(t^{-1}\) at infinity. Apply
Lemma~\ref{lem:regularization} to the probabilities \((w_i)\)
and the exponents \(r,q\). The resulting probabilities \((v_i)\)
and fixed diagonal map \(V\) satisfy
\begin{align*}
 \|VU\widehat A x\|_{q;v}
 &\leq(1+\log(2k))^{1/q}
       \|VU\widehat A x\|_{q,\infty;v}\\
 &\leq C(r,q)(1+\log(2k))^{1/q}
       \|U\widehat A x\|_{q,\infty;w}\\
 &\leq C(r,q)C_1(r,q)DT(1+\log(2m))^{1/q}\|x\|_X.
\end{align*}
On the other hand, \(v\) is a probability and \(r<q\), so
\[
 \|VU\widehat A x\|_{q;v}
 \geq\|VU\widehat A x\|_{r;v}
 =\|U\widehat A x\|_{r;w}
 =\|Ax\|_r\geq\|x\|_X.
\]
For \(y\in A(X)\), define
\[
 Ry=\bigl(v_i^{1/q}[VU(m^{1/r}y)](i)\bigr)_{i=1}^k
          \oplus0\in\ell_q^m.
\]
This map is \(\K\)-linear. Substituting \(y=Ax\) gives the
required bounds for \(RA\), and also proves that \(R\)
is injective on \(A(X)\). The density,
regularized probabilities, and coordinate multipliers are fixed
for the whole space, and at most \(m\) coordinates are used.
This proves (ii).
\end{proof}

\section{Proof of the main theorem}\label{sec:mainproof}
The same map into a lower-exponent coordinate space will serve for
both assertions. Its distortion is independent of the improved
type parameters. Those enter only when the exponent is restored.

\begin{proof}[Proof of Theorem~\ref{thm:targets}]
We first assume \(1<\mu\leq2\), and put
\[
 a=p/2,\qquad r=3p/8,\qquad m=\lfloor\mu d\rfloor.
\]
Lemma~\ref{lem:preliminary} gives a \(\K\)-linear embedding
\(J:X\to\ell_a^N\), which we normalize so that
\[
 \|x\|_p\leq\|Jx\|_a\leq D_a\|x\|_p,
 \qquad N\leq C_0(p)d.
\]
Since \(a/2<r<a\), Proposition~\ref{prop:compression} applies
to \(E=JX\). After normalization its embedding
\(S:E\to\ell_r^m\) satisfies
\[
 \|z\|_a\leq\|Sz\|_r\leq D_c\|z\|_a,
 \qquad D_c=K(a,r,C_0(p),\mu).
\]
Set \(A=SJ:X\to\ell_r^m\) and \(Y=A(X)\), with the
quasi-norm inherited from \(\ell_r^m\). Then
\begin{equation}\label{eq:main-composite}
 \|x\|_p\leq\|Ax\|_r\leq D_0\|x\|_p,
 \qquad D_0=D_aD_c.
\end{equation}
By Remark~\ref{rem:excess}, increasing a constant depending only
on \(p\) gives
\begin{equation}\label{eq:main-lower-distortion}
 D_0\leq\exp\!\left(\frac{C_0'(p)}{\mu-1}\right).
\end{equation}
Indeed, \(a,r,C_0(p)\) depend only on \(p\), and the fixed
factor \(D_a\) can be absorbed because \(\mu-1\leq1\).

Suppose now that \(T_q(X)\leq T\), where \(p<q\leq2\).
Lemma~\ref{lem:restore}(i) applies directly to \(X\) and \(A\),
with \(D=D_0\), source exponent \(r\), target exponent \(p\),
and type exponent \(q\). Its hypotheses hold because
\(r<p<q\leq2\) and \eqref{eq:main-composite} gives the required
upper and lower bounds. It yields a \(\K\)-linear embedding
\(R:Y\to\ell_p^m\) such that
\[
 \|x\|_p\leq\|RAx\|_p\leq C(r,p,q)D_0T\|x\|_p.
\]
Consequently, \(RA:X\to\ell_p^m\) satisfies
\[
 \operatorname{dist}(RA)\leq C(r,p,q)D_0T
 \leq C(p,q)T\exp\!\left(\frac{C_0'(p)}{\mu-1}\right).
\]
This proves (a), including its quantitative bound. The same argument applies when
\(q=2\), since the preliminary reduction used only the automatic
type \(p\) of \(X\).

For (b), the \(p\)-triangle inequality gives, for every finite
sequence \((x_i)\subset X\) and every choice of signs,
\[
 \left\|\sum_i\varepsilon_i x_i\right\|_p
 \leq\left(\sum_i\|x_i\|_p^p\right)^{1/p}.
\]
Taking the second moment shows \(T_p(X)\leq1\).
Apply Lemma~\ref{lem:restore}(ii) to \(X\) and \(A\), with
\(D=D_0\), type exponent equal to \(p\), and type bound one.
Its source exponent is \(r<p<1\). Denote the resulting
restoration map by \(R:Y\to\ell_p^m\). Its bounds give
\[
 \|x\|_p\leq\|RAx\|_p
 \leq C(p)D_0(1+\log(2m))^{1/p}\|x\|_p,
\]
and hence
\[
 \operatorname{dist}(RA)
 \leq C(p)D_0(1+\log(2m))^{1/p}.
\]
Since \(m\leq2d\) and \(d\geq1\),
\[
 1+\log(2m)\leq1+\log(4d)
 \leq\left(1+\frac{1+\log4}{\log3}\right)\log(d+2).
\]
Combine this with \eqref{eq:main-lower-distortion} and absorb
the remaining factor depending only on \(p\) into the exponential.
After increasing \(C(p)\), we obtain the second bound in
\eqref{eq:main-quantitative}.

Finally, if \(\mu>2\), apply the result with \(\mu=2\) and
append zeros to reach \(\lfloor\mu d\rfloor\) coordinates.
This does not change the distortion. All maps in the construction
are \(\K\)-linear. Proposition~\ref{prop:compression} includes
the small-dimensional cases, so the argument covers every
\(d\geq1\).
\end{proof}

\begin{remark}\label{rem:excess-lower}
Some divergence of the distortion as \(\mu\downarrow1\)
is unavoidable even for Hilbertian sources. To make this precise,
let \(F_p(\mu)\) be any bound that works in
Theorem~\ref{thm:targets}(a) for every \(d\), with \(q=2\)
and \(T=1\). Put \(\varepsilon=\mu-1\), and assume
\(0<\varepsilon\leq1/4\). Choose
\[
 d=\left\lfloor\frac1{2\varepsilon}\right\rfloor.
\]
Then \(d\geq1/(4\varepsilon)\) and
\(0<\varepsilon d\leq1/2\), so
\(\lfloor\mu d\rfloor=d\).

The space \(\ell_2^d(\K)\) is linearly isometric to a subspace
of \(L_p(\K)\) and has type \(2\) constant one. Indeed, let
\(G_1,\ldots,G_d\) be independent standard real Gaussians when
\(\K=\R\), and independent circular complex Gaussians when
\(\K=\mathbb C\). With
\(\gamma_p=(\E|G_1|^p)^{1/p}\), the map
\[
 x\longmapsto\gamma_p^{-1}\sum_{i=1}^d x_iG_i
\]
is \(\K\)-linear and isometric, since its unnormalized image
has the distribution of \(\|x\|_2G_1\).

Any linear embedding \(B:\ell_2^d\to\ell_p^d\) is onto.
Testing the type \(2\) inequality on the standard coordinate
vectors of \(\ell_p^d\) gives
\[
 T_2(\ell_p^d)\geq
 \frac{\bigl(\E\|\sum_{i=1}^d\varepsilon_i e_i\|_p^2\bigr)^{1/2}}
      {(\sum_{i=1}^d\|e_i\|_p^2)^{1/2}}
 =d^{1/p-1/2}.
\]
The transport inequality \eqref{eq:type-transport}, applied to
\(B\), therefore implies
\[
 \operatorname{dist}(B)\geq d^{1/p-1/2},
 \qquad
 F_p(1+\varepsilon)\geq
 (4\varepsilon)^{-(1/p-1/2)}.
\]
Thus a polynomial divergence is necessary for a bound uniform in
all dimensions, whereas our upper bound is exponential.
This example uses the dimensions for which rounding leaves
exactly \(d\) target coordinates. It does not give the same
obstruction in the regime \(\varepsilon d\to\infty\), nor does
it establish optimality of our upper bound.
\end{remark}

\begin{remark}\label{rem:logarithm}
The logarithm in Theorem~\ref{thm:targets}(b) enters only when a
weak estimate at the available type exponent is converted into a
strong estimate. The compression and the preliminary reduction
have distortion independent of \(d\). After regularization, the
endpoint tail integral is truncated at
\((2k)^{1/q}\|Vf\|_{q,\infty;v}\), and integrating \(t^{-1}\)
gives \(\log(2k)\). This explains the bound obtained by the method.
It does not establish that the logarithmic factor is necessary for
arbitrary subspaces of \(L_p\). We prove no matching lower bound.
For example, the particular source \(\ell_p^d\) embeds isometrically
into \(\ell_p^m\) for every \(m\geq d\) by appending zeros.
The issue of improving the bound concerns a uniform estimate over
all \(d\)-dimensional subspaces, not such individual examples.
\end{remark}

\end{document}